\documentclass[graybox]{svmult}

\usepackage[all,cmtip]{xy}

\usepackage{mathptmx}       
\usepackage{helvet}         
\usepackage{courier}        
 \usepackage{type1cm}        
\usepackage{makeidx}         
\usepackage{graphicx}        
\usepackage{multicol}        
\usepackage[bottom]{footmisc}
\usepackage{amsmath,amssymb,bm} 

\makeindex             

\usepackage{amsmath,amssymb,bm}

\newcommand{\cK}{{\mathcal{K}}}

\newcommand{\C}{{\mathbb{C}}} 
\newcommand{\R}{{\mathbb{R}}} 

\DeclareMathOperator{\Val}{{\rm Val}}

\DeclareMathOperator{\Curv}{Curv}

\begin{document}

\title*{Valuations and curvature measures on complex spaces}
\author{Andreas Bernig}
\institute{Andreas Bernig \at Institut f\"ur Mathematik, Goethe-Universit\"at Frankfurt,
Robert-Mayer-Str. 10, 60054 Frankfurt, Germany, \email{bernig@math.uni-frankfurt.de\\Supported by DFG grant BE 2848/5-2.}
}
%
%
\maketitle

\abstract{We survey recent results in hermitian integral geometry, i.e. integral geometry on complex vector spaces and complex space forms. We study valuations and curvature measures on complex space forms and describe how the global and local kinematic formulas on such spaces were recently obtained. While the local and global kinematic formulas in the Euclidean case are formally identical, the local formulas in the hermitian case contain strictly more information than the global ones. Even if one is only interested in the flat hermitian case, i.e. $\C^n$, it is necessary to study the family of all complex space forms, indexed by the holomorphic curvature $4\lambda$, and the dependence of the formulas on the parameter $\lambda$. We will also describe Wannerer's recent proof of local additive kinematic formulas for unitarily invariant area measures.}

\section{Introduction}

Hermitian integral geometry is a relatively old subject, with early contributions by Blaschke \cite{blaschke39b}, Rohde \cite{rohde40}, Santal\'o \cite{santalo52}, Shifrin \cite{shifrin81}, Gray \cite{gray_book}, Griffiths \cite{griffiths78} and others, and more recent works by Tasaki \cite{tasaki00, tasaki03}, Park \cite{park02}, and Abardia-Gallego-Solanes \cite{abardia_gallego_solanes}.

A systematic study of this subject was completed only recently. This is partly due to the fact that Alesker's algebraic theory of translation invariant valuations as well as his theory of valuations on manifolds, both being indispensable tools in this line research, were not available before. The main recent results of Hermitian integral geometry are the explicit description of global kinematic formulas in hermitian spaces \cite{bernig_fu_hig}, local and global kinematic formulas on complex space forms 
\cite{bernig_fu_solanes} and kinematic formulas for unitarily invariant area measures \cite{wannerer_area_measures}. 

These three papers cover almost 200 pages. The aim of the present paper is to state most of the main theorems and to give some ideas about their proofs. More generally, we want to give a taste of this beautiful theory, which involves tools and results from convex geometry, differential geometry, geometric measure theory, representation theory and algebraic geometry. 

In Section \ref{sec_translation_inv} we sketch some fundamental facts and definitions of Alesker's theory of translation invariant valuations, in particular we introduce product, convolution and Alesker-Fourier transform. Section \ref{sec_curv_meas} deals with valuations on manifolds; and introduces the important notion of a smooth curvature measure on a manifold. In Section \ref{sec_kinforms} we introduce the general framework to study local and global kinematic formulas on space forms and state the important transfer principle. 

The remaining sections deal with the hermitian case. Section \ref{sec_hermitian_flat} is a survey on the results from \cite{bernig_fu_hig}, the most important theorem in this section being Theorem \ref{thm_pkf_flat}. The curved case, i.e. the case of complex space forms, is treated in Section \ref{sec_hermitian_curved}. The main results are Theorems \ref{thm_isomorphism} and  \ref{thm_plkf}. Finally, the theory of unitarily invariant area measures is sketched in Section \ref{sec_area_meas}. This theory has some similar features as the theory of curvature measures. We do not explicitly state Wannerer's theorem (the local additive kinematic formulas for unitarily invariant area measures), but sketch the way they are obtained.

\section{Translation invariant valuations on vector spaces}
\label{sec_translation_inv}

In this section, we will introduce the basic notation and state some fundamental theorems, some of which will be explained in more detail in \cite{alesker_lnm}. 

\subsection{McMullen's decomposition and Hadwiger's theorem}

Let $\cK^n$ be the space of compact convex bodies in $\R^n$. A (real-valued) valuation on $\R^n$ is a functional $\mu:\cK^n \to \R$ which satisfies
\begin{displaymath}
 \mu(K \cup L)+\mu(K \cap L)=\mu(K)+\mu(L)
\end{displaymath}
whenever $K,L,K \cup L \in \cK^n$. 

Continuity of valuations is with respect to the Hausdorff topology. A valuation is said to be translation invariant, if $\mu(K+x)=\mu(K)$ for all $K \in \cK^n, x \in \R^n$. The vector space of all continuous translation invariant valuations is denoted by $\Val$.  

 If $G$ is some group acting linearly on $\R^n$, a valuation $\mu$ is called $G$-invariant if $\mu(gK)=\mu(K)$ for all $K \in \cK^n, g \in G$.

A valuation $\mu$ is called homogeneous of degree $k$ if $\mu(tK)=t^k \mu(K), t>0$. It is even/odd if $\mu(-K)=\mu(K)$ (resp. $\mu(-K)=-\mu(K)$) for all $K$. The corresponding space is denoted by $\Val_k^\pm$.

A fundamental theorem concerning translation invariant valuations is McMullen's decomposition. 
\begin{theorem}[McMullen, \cite{mcmullen77}]
 \begin{displaymath}
  \Val=\bigoplus_{\substack{k=0,\ldots,n\\ \epsilon=\pm}} \Val_k^\epsilon.
 \end{displaymath}
\end{theorem}

The space $\Val_n$ is one-dimensional and spanned by the volume, while $\Val_0$ is spanned by the Euler characteristic $\chi$.  

The space of $\mathrm{SO}(n)$-invariant valuations was described by Hadwiger. Let $\mu_k$ denote the $k$-th intrinsic volume, see \cite{klain_rota, schneider_book14}. 

\begin{theorem}[Hadwiger, \cite{hadwiger_vorlesung}, see also \cite{klain95, klain_rota}]
The intrinsic volumes $\mu_0,\ldots,\mu_n$ form a basis of $\Val^{\mathrm{SO}(n)}$.
\end{theorem}

Note that without any $G$-invariance, the spaces $\Val_k^\pm$ are infinite-dimensional (except for $k=0,n$). Using McMullen's theorem, one can show that they admit a Banach space structure. 

Even valuations are easier to understand than odd ones thanks to Klain's embedding result \cite{klain00}. If $\phi \in \Val_k^+$, then the restriction of $\phi$ to a $k$-dimensional plane $E \in \mathrm{Gr}_k$ is a multiple $\mathrm{Kl}_\phi(E)$ of the Lebesgue measure. Klain proved that the map $\Val_k^+ \to C(\mathrm{Gr}_k), \phi \mapsto \mathrm{Kl}_\phi$ is injective. 

\subsection{Alesker's irreducibility theorem, product and convolution of valuations}

The group $\mathrm{GL}(n)$ acts on $\Val$ by 
\begin{displaymath}
 (g\mu)(K):=\mu(g^{-1}K).
\end{displaymath}
Obviously, degree and parity are preserved under this action. 

One of the most fundamental and influential theorems of modern integral geometry is the following. 

\begin{theorem}[Alesker, \cite{alesker_mcullenconj01}]
 The spaces $\Val_k^\epsilon, k=0,\ldots,n, \epsilon=\pm$ are irreducible $\mathrm{GL}(n)$-representations. 
\end{theorem}

Since we are in an infinite-dimensional situation, this means that every non-trivial $\mathrm{GL}(n)$-invariant subspace of $\Val_k^\epsilon$ is dense. 

\begin{corollary}
 Valuations of the form $K \mapsto \mathrm{vol}(K+A)$ with $A \in \cK^n$ span a dense subspace inside $\Val$.
\end{corollary}
  
\begin{definition}
 A valuation $\mu \in \Val$ is called smooth, if the map 
\begin{align*}
 \mathrm{GL}(n) & \to \Val\\
g & \mapsto g\mu
\end{align*}
is smooth as a map from the Lie group  $\mathrm{GL}(n)$ to the Banach space $\Val$.
\end{definition}

Smooth valuations form a dense subspace $\Val^\infty \subset \Val$, which has a natural Fr\'echet space structure. 

An example of a smooth valuation is a valuation of the form $K \mapsto \mathrm{vol}(K+A)$, where $A$ is a convex body with smooth boundary and positive curvature. 

Based on Alesker's Irreducibility theorem, a rich algebraic structure was introduced on $\Val^\infty$ in recent years. 
\begin{theorem}[Alesker \cite{alesker03_un, alesker04_product, alesker_fourier}, Bernig-Fu \cite{bernig_fu06}]
\label{thm_alg_structures}
\begin{enumerate}
\item[(Alesker)] There exists a unique continuous bilinear product on $\Val^\infty$ such that if 
\begin{displaymath}
 \mu_i(K)=\mathrm{vol}_n(K+A_i), i=1,2,
\end{displaymath}
with convex bodies $A_i$ with smooth boundary and positive curvature, then 
\begin{displaymath}
 \mu_1 \cdot \mu_2(K)=\mathrm{vol}_{2n}(\Delta K+A_1 \times A_2),
\end{displaymath}
where $\Delta:\R^n \to \R^n \times \R^n$ is the diagonal embedding. 
\item[(Bernig-Fu)]  There exists a unique continuous bilinear product on $\Val^\infty$ such that if 
\begin{displaymath}
 \mu_i(K)=\mathrm{vol}_n(K+A_i), i=1,2,
\end{displaymath}
with convex bodies $A_i$ with smooth boundary and positive curvature, then 
\begin{displaymath}
 \mu_1 * \mu_2(K)=\mathrm{vol}_{n}(K+A_1 + A_2).
\end{displaymath}
\item[(Alesker)] There exists a Fourier-type transform $\mathbb F:\Val^\infty \to \Val^\infty$ (called Alesker-Fourier transform) such that 
\begin{displaymath}
 \mathbb F(\mu_1 \cdot \mu_2) = \mathbb F \mu_1 * \mathbb F \mu_2, \quad \mu_1,\mu_2 \in \Val^\infty
\end{displaymath}
and 
\begin{displaymath}
 \mathbb F^2=\epsilon \quad \text{ on } \Val_k^\epsilon.
\end{displaymath}
\end{enumerate}
\end{theorem}

We note that, up to some minor modifications, these algebraic structures do not depend on the choice of the Euclidean scalar product. We refer to \cite{alesker_lnm} for more details.

In Section \ref{sec_kinforms} it will become transparent that product and convolution are very closely related to kinematic formulas. 

\begin{theorem}[\cite{alesker04_product}]
 The map
 \begin{displaymath}
  \Val^\infty \otimes \Val^\infty \to \Val_n \cong \mathbb{R}, \phi_1 \otimes \phi_2 \mapsto (\phi_1 \cdot \phi_2)_n
 \end{displaymath}
is perfect and thus induces an injective map with dense image
\begin{displaymath}
 \mathrm{PD}: \Val^\infty \to (\Val^\infty)^*
\end{displaymath}

\end{theorem}

A smooth valuation can be described by a pair of differential forms as follows. Let $K$ be a convex body and $N(K)$ its normal cycle \cite{zaehle86}. As a set, it consists of all pairs $(x,v)$, where $x \in \partial K$ and $v$ is an outer normal vector of $K$ at $x$. This set is an $(n-1)$-dimensional Lipschitz submanifold of the sphere bundle $S\R^n=\R^n \times S^{n-1}$, which can be endowed in a canonical way with an orientation. Given a smooth $(n-1)$-form $\omega$ on $S\R^n$, we can integrate it over $N(K)$. 
\begin{proposition} \label{prop_smooth_vals}
 A valuation $\mu \in \Val$ is smooth if and only if there are translation invariant differential forms $\omega \in \Omega^{n-1}(S\R^n), \phi \in \Omega^n(\R^n)$ such that 
 \begin{displaymath}
  \mu(K)=\int_{N(K)} \omega+\int_K \phi
 \end{displaymath}
for all $K$. 
\end{proposition}
These forms are not unique, see Theorem \ref{thm_kernel_thm} below. 

\section{Valuations and curvature measures on manifolds}
\label{sec_curv_meas}

In this section, we define smooth valuations on manifolds and describe them in terms of differential forms. Then we extend the product from the previous section to such valuations. Smooth curvature measures on valuations are introduced as measure-valued smooth valuations. They form a \emph{module} over the space of smooth valuations. We refer to \cite{ alesker_val_man1, alesker_val_man2, alesker_survey07, alesker_val_man4, alesker_intgeo, alesker_barcelona, alesker_bernig, alesker_val_man3,   bernig_broecker07} for valuations on manifolds and to \cite{bernig_fu_solanes,park02} for curvature measures on manifolds. 
 
\subsection{Smooth valuations on manifolds, the Rumin operator and the product structure}

If we want to define the notion of valuation on a smooth manifold of dimension $n$, an obvious obstacle is that the notion of convex set is not available. In the presence of a Riemannian metric, one can define convex sets, but their behaviour is too wild to be a good substitute of the notion of convexity on Euclidean space. Instead, we are led to define valuations on some other class of reasonable sets. One possibility is to use sets of positive reach or variants of it. Another possibility is to use compact submanifolds with corners, also called simple differentiable polyhedra. Each simple differential polyhedron has a conormal cycle, which is an $(n-1)$-dimensional Lipschitz manifold in the cosphere bundle $S^*M=\{(x,[\xi]): x \in M, \xi \in T^*_xM \setminus \{0\}\}$, where the equivalence classes are taken with respect to the relation $\xi_1 \sim \xi_2$ if and only if $\xi_2 =  \lambda \xi_1$ for some $\lambda>0$. 

Here we will not pay much attention to the precise class of sets, since the kinematic formulas in the different settings are formally identical. We refer to \cite{fu_lnm} for a more thorough study of this question. 

The second difference with the flat case is that we do not really define valuations, but only smooth valuations. The definition is inspired by Proposition \ref{prop_smooth_vals}. 

\begin{definition}
 A smooth valuation on an $n$-dimensional manifold $M$ is a functional $\mu$ on the space $\mathcal{P}(M)$ of simple differentiable polyhedra which has the form 
 \begin{displaymath}
  \mu(P)=\int_{N(P)} \omega+\int_P \phi
 \end{displaymath}
with smooth forms $\omega \in \Omega^{n-1}(S^*M), \phi \in \Omega^n(M)$. 
Here $N(P) \subset S^*M$ denotes the conormal cycle of $P$. The space of smooth valuations on $M$ is denoted by $\mathcal{V}(M)$. 
\end{definition}

Note that no invariance is assumed in this definition. 

Taking $\omega=0$ and $\phi$ a volume form (say with respect to a Riemannian metric), we see that the Riemannian volume is a smooth valuation. The Euler characteristic is another example. In fact, Chern \cite{chern44,chern45}, when proving the famous Chern-Gauss-Bonnet theorem, constructed a pair $(\omega,\phi)$ as above. By Proposition \ref{prop_smooth_vals}, each smooth translation invariant valuation on the vector space $\R^n$ can also be considered as a smooth valuation on the manifold $\R^n$.

An important point to note is that the pair of forms $(\omega,\phi)$ is not unique. Since the conormal cycle is closed, it annihilates exact differential forms. Moreover, it is Legendrian, i.e. it annihilates forms which vanish on the contact distribution in $S^*M$. The kernel of the map which associates to a pair of forms the corresponding smooth valuation is given in the following theorem.
\begin{theorem}[\cite{bernig_broecker07}] \label{thm_kernel_thm}
 A pair of forms $(\omega,\phi)$ induces the trivial valuation if and only if 
 \begin{enumerate}
  \item $D\omega+\pi^*\phi=0$,
  \item $\pi_* \omega=0$.
 \end{enumerate}
Here $D:\Omega^{n-1}(S^*M) \to \Omega^n(S^*M)$ is a certain second order differential operator, called Rumin operator \cite{rumin94}, $\pi^*$ denotes pull-back and $\pi_*$ push-forward (or fiber integration) with respect to the projection map $\pi:S^*M \to M$. 
\end{theorem}

Any operation on pairs of forms $(\omega,\phi)$ which is compatible with the kernel described in Theorem \ref{thm_kernel_thm} thus induces an operation on smooth valuations. An easy example is given by the Euler-Verdier involution, which on the level of forms is given by $(\omega,\phi) \mapsto ((-1)^n a^*\omega,(-1)^n \phi)$, where $a:S^*M \to S^*M, (x,[\xi]) \mapsto (x,[-\xi])$ is the natural involution (anti-podal map). 

A much more involved example is the product of smooth valuations on manifolds. The complete formula is rather technical and uses certain blow-up spaces whose definition we prefer to omit. 
\begin{theorem}[\cite{alesker_bernig, alesker_val_man3}]
\label{thm_product_formula}
 The space of smooth valuations on a manifold admits a product structure which, on the level of forms, is given by a formula of the type 
 \begin{displaymath}
  (\omega_1,\phi_1) \cdot (\omega_2,\phi_2) \mapsto (Q_1(\omega_1,\omega_2),Q_2(\omega_1,\phi_1,\omega_2,\phi_2))
 \end{displaymath}
Here $Q_1,Q_2$ are certain explicitly known operators on differential forms involving some Gelfand transform on a blow-up space. This product is commutative, associative, has thje Euler characteristic as unit, and is compatible with restrictions to submanifolds. 
\end{theorem}

A nice interpretation of this formula was recently given in \cite{fu_alesker_product}. Of course, in the (very) special case of smooth translation invariant valuations on $\R^n$, the new product coincides with the product from Theorem \ref{thm_alg_structures}. 

Let us now describe a version of Poincar\'e duality on the level of valuations on manifolds, which was introduced in \cite{alesker_val_man4}. Given any compactly supported smooth valuation $\mu$ on a manifold $M$, we may evaluate it at the manifold $M$ to obtain a real number denoted by $\int \mu$. 

If we denote by $\mathcal{V}_c(M)$ the compactly supported smooth valuations on $M$, then we get a pairing 
\begin{equation} \label{eq_poincare_compact}
 \mathcal{V}(M) \times \mathcal{V}_c(M) \to \R, (\mu_1,\mu_2) \mapsto \int \mu_1 \cdot \mu_2.
\end{equation}
Alesker has shown that this pairing is perfect, i.e. it induces an injective map $\mathcal{V}(M) \to \mathcal{V}_c(M)^*$ with dense image. This fact is important in connection with \emph{generalized valuations on manifolds} (see \cite{alesker_val_man4, alesker_bernig, alesker_faifman}), but also for kinematic formulas \cite{alesker_bernig, bernig_fu_solanes}.

\subsection{Curvature measures, module structure}

Roughly speaking, a curvature measure is a valuation with values in the space of signed measures. An example are Federer's curvature measures $\Phi_0,\ldots,\Phi_n$ in $\R^n$, which (up to scaling) are localizations of the intrinsic volumes. 

\begin{definition}[\cite{bernig_fu_solanes}]
 A smooth curvature measure on an $n$-dimensional manifold $M$ is a functional $\Phi$ of the form 
\begin{displaymath}
 \Phi(P,B)=\int_{N(P) \cap \pi^{-1}B} \omega+\int_{P \cap B} \phi,
\end{displaymath}
where $P$ is a simple differentiable polyhedron, $B \subset M$ a Borel subset, and $\omega \in \Omega^{n-1}(S^*M), \phi \in \Omega^n(M)$. The globalization of $\Phi$ is the smooth valuation $\mathrm{glob}\ \Phi(P):=\Phi(P,M)$. The space of all smooth curvature measures on $M$ is denoted by $\mathcal{C}(M)$.
\end{definition}

Let us explain where the name ``curvature measure`` comes from. For this, suppose that $P$ is a compact smooth submanifold of a Riemannian manifold $M$. The geometry of the second fundamental form of $P$ is then encoded in the conormal cycle $N(P)$. If $\Phi$ is a smooth curvature measure, then the measure $\Phi(P,\cdot)$ is obtained by integration of some polynomial function on the second fundamental form, hence by some curvature expression. 

As in Theorem \ref{thm_kernel_thm}, one may describe the kernel of the map which associates to a pair of differential forms a curvature measure. The pair $(\omega,\phi)$ induces the trivial curvature measure if and only if $\phi=0$ and $\omega$ belongs to the ideal generated by $\alpha$ and $d\alpha$, where $\alpha$ is the contact form on $S^*M$.

It follows that the map $\mathrm{glob}:\mathcal{C}(M) \to \mathcal{V}(M)$ is surjective, but not injective. For example, if $\omega$ is an exact form which is not contained in the ideal generated by $\alpha$ and $d\alpha$, then $(\omega,0)$ defines a non-zero curvature measure whose globalization vanishes.

A more general globalization map is obtained as follows. If $f$ is a smooth function on $M$, we may define a smooth valuation $\mathrm{glob}_f \Phi$ by 
\begin{displaymath}
 \mathrm{glob}_f \Phi(P):=\int_{N(P)} \pi^*f \ \omega + \int_P f \phi. 
\end{displaymath}
In the particular case $f \equiv 1$, this is just the map $\mathrm{glob}$. 

While we can multiply smooth valuations on a manifold, there seems to be no reasonable product structure on the space of smooth curvature measures. However, curvature measures form a module over valuations. 

\begin{theorem}[\cite{bernig_fu_solanes}, based on \cite{alesker_bernig}]
\label{thm_module}
The space $\mathcal{C}(M)$ is a module over the algebra $\mathcal{V}(M)$. More precisely, given $\mu \in \mathcal{V}(M), \Phi \in \mathcal{C}(M)$, there exists a unique curvature measure $\mu \cdot \Phi \in \mathcal{C}(M)$ such that $\mathrm{glob}_f(\mu \cdot \Phi)=\mu \cdot \mathrm{glob}_f \Phi$ for all smooth functions $f$ on $M$. 
\end{theorem}

The proof follows rather easily from Theorem \ref{thm_product_formula}. 

\section{Global and local kinematic formulas, the transfer principle}
\label{sec_kinforms}

In this section, we study the global and local kinematic operators on isotropic spaces, their link to the algebra and module structure from the previous section, and the important transfer principle relating local kinematic formulas on flat and curved spaces. 
 
\subsection{Fundamental theorem of algebraic integral geometry I: flat case}

\begin{definition}
 The space of smooth and translation invariant curvature measures on $\mathbb R^n$ is denoted by $\mathrm{Curv}$. If $G$ acts linearly on $\R^n$, then $\mathrm{Curv}^G$ is the subspace of $G$-invariant elements.
\end{definition}

As an example, $\mathrm{Curv}^{\mathrm{SO}(n)}$ has Federer's curvature measures $\Phi_0,\ldots,\Phi_n$ as basis. In particular, the restricted globalization map $\mathrm{glob}:\mathrm{Curv}^{\mathrm{SO}(n)} \to \mathrm{Val}^{\mathrm{SO}(n)}$ is a bijection.

\begin{theorem}[\cite{alesker_survey07}]
Let $G$ be a subgroup of $\mathrm{O}(n)$. Then $\Val^G$ is finite-dimensional if and only if $G$ acts transitively on the unit sphere. In this case, $\Val^G \subset \Val^\infty$ and  $\mathrm{Curv}^G$ is finite-dimensional as well.   
\end{theorem}

Let $G$ be such a group and denote by $\bar G$ the group generated by $G$ and translations (endowed with a convenient invariant measure). Let $\phi_1,\ldots,\phi_N$ be a basis of $\Val^G$. Then there are constants $c_{kl}^i$ such that the kinematic formulas
\begin{displaymath}
 \int_{\bar G} \phi_i(K \cap \bar g L) d\bar g=\sum_{k,l} c_{k,l}^i \phi_k(K) \phi_l(L), \quad K,L \in \cK^n
\end{displaymath}
hold. 

Similarly, if $\Phi_1,\ldots,\Phi_M$ is a basis of $\mathrm{Curv}^G$, there are constants $d_{k,l}^i$ such that 
\begin{displaymath}
 \int_{\bar G} \Phi_i(K \cap \bar gL,B_1 \cap \bar g B_2) d\bar g=\sum_{k,l} d_{k,l}^i \Phi_k(K,B_1)\Phi_l(L,B_2)
\end{displaymath}
holds for all convex bodies $K,L$ and all Borel subsets $B_1,B_2 \subset \R^n$. The proof is contained in \cite{fu90}.

We call the corresponding operators 
\begin{align*}
 k_G: \Val^G & \to \Val^G \otimes \Val^G \\
\phi_i & \mapsto \sum_{k,l} c_{k,l}^i \phi_k \otimes \phi_l
\end{align*}
and 
\begin{align*}
 K_G: \mathrm{Curv}^G & \to \mathrm{Curv}^G \otimes \mathrm{Curv}^G \\
\Phi_i & \mapsto \sum_{k,l} d_{k,l}^i \Phi_k \otimes \Phi_l
\end{align*}
the global and local kinematic operators. They are independent of the choice of the bases. Moreover, Fubini's theorem and the unimodularity of $G$ imply that $k_G, K_G$ are cocommutative coassociative coproducts. 

Since the diagram 
\begin{displaymath}
 \xymatrix{\mathrm{Curv}^G \ar[r]^-{K_G} \ar[d]^{\mathrm{glob}} & \mathrm{Curv}^G \otimes \mathrm{Curv}^G \ar[d]^{\mathrm{glob} \otimes \mathrm{glob}} \\ 
\Val^G \ar[r]^-{k_G} & \Val^G \otimes \Val^G}
\end{displaymath}
clearly commutes and $\mathrm{glob}$ is surjective, but usually non-injective, the local operator kinematic operator contains more information than the global one. 

There is another kinematic operator which lies between these two, called semi-local kinematic operator, which is defined by 
\begin{displaymath}
 \bar k_G:=(\mathrm{id} \otimes \mathrm{glob}) \circ K_G:\mathrm{Curv}^G \to \mathrm{Curv}^G \otimes \mathrm{Val}^G.
\end{displaymath}

The following theorem, although rather easy to prove, is fundamental for our algebraic understanding of kinematic formulas. 
\begin{theorem}[\cite{fu06, bernig_fu06}]
\label{thm_ftaig}
 Let $m:\mathrm{Val}^G \otimes \mathrm{Val}^G \to \mathrm{Val}^G$ denote the restriction of the product to $G$-invariant valuations and $m^*:\mathrm{Val}^{G*} \to \mathrm{Val}^{G*} \otimes \mathrm{Val}^{G*}$ its adjoint. Then the following diagram commutes
 \begin{equation} \label{eq_ftaig}
  \xymatrix{\mathrm{Val}^G \ar[r]^-{k_G} \ar[d]^{\mathrm{PD}} & \mathrm{Val}^G \otimes \mathrm{Val}^G \ar[d]^{\mathrm{PD} \otimes \mathrm{PD}}\\
  \mathrm{Val}^{G*} \ar[r]^-{m^*} & \mathrm{Val}^{G*} \otimes \mathrm{Val}^{G*}}
 \end{equation}
\end{theorem}

Thus knowing the product structure, we can (at least in principle) write down the global kinematic formulas. 

\subsection{The transfer principle}

\begin{definition}
An isotropic space is a pair $(M,G)$, where $M$ is a Riemannian manifold and $G$ is a Lie subgroup of the isometry group of $M$ which acts transitively on the sphere bundle $SM$.  
\end{definition}

\begin{theorem}[\cite{alesker_survey07, bernig_fu_solanes}]
Let $(M,G)$ be an isotropic space. Then $\mathcal{V}(M)^G$ and $\mathcal{C}(M)^G$ are finite-dimensional.   
\end{theorem}

If $(M,G)$ is an isotropic space, there exist local kinematic formulas as follows. Let $X,Y$ be sufficiently nice sets (compact submanifolds with corners will be enough for our purpose) and let $B_1,B_2 \subset M$ be Borel subsets. If $\Phi_1,\ldots,\Phi_N$ is a basis of $\mathcal{C}(M)^G$, then there are constants $d_{k,l}^i$ such that 
\begin{displaymath}
\int_G \Phi_i(X \cap gY,B_1 \cap gB_2) dg=\sum_{k,l} d_{k,l}^i \Phi_k(X,B_1)\Phi_l(Y,B_2). 
\end{displaymath}
The existence of such formulas is harder to prove than in the flat case (compare \cite{bernig_fu_solanes, fu90}). One reason for this is that it is unknown whether there are $G$-invariant continuous, but non-smooth valuations (even in the simplest case of the sphere). 

As in the flat case, we obtain a map $K_G:\mathcal{C}(M)^G \to \mathcal{C}(M)^G \otimes \mathcal{C}(M)^G$. Semi-local and global kinematic formulas are defined similarly, the corresponding operators are denoted by $\bar k_G,k_G$.  Then $(\mathcal{V}(M)^G,k_G)$ and $(\mathcal{C}(M)^G,K_G)$ are cocommutative coassociative coalgebras and the following diagram is commutative.
\begin{equation} \label{eq_local_and_global_formulas_curved}
\xymatrix{ \mathcal{C}(M)^G \ar[r]^-{K_G} \ar[d]^{\mathrm{id}} &  \mathcal{C}(M)^G \otimes \mathcal{C}(M)^G  \ar[d]^{\mathrm{id} \otimes \mathrm{glob}}\\
\mathcal{C}(M)^G \ar[r]^-{\bar k_G} \ar[d]^{\mathrm{glob}} &  \mathcal{C}(M)^G \otimes \mathcal{V}(M)^G  \ar[d]^{\mathrm{glob} \otimes \mathrm{id}}\\
\mathcal{V}(M)^G \ar[r]^-{k_G} & \mathcal{V}(M)^G \otimes \mathcal{V}(M)^G}
\end{equation}

The transfer principle below will apply to the following families of isotropic spaces.
\begin{equation} \label{table_isotropic_spaces}
 \begin{array}{l | l | l | l}
\text{division algebra} & \text{positively curved} & \text{flat} & \text{negatively curved} \\  \hline
  \R & (S^n,\mathrm{SO(n+1)}) & (\R^n,\overline{\mathrm{SO}(n)} & (\mathbb{H}^n,\mathrm{SO}(n,1))\\
 \C & (\mathbb{C}P^n,\mathrm{PU}(n+1)) & (\C^n,\overline{\mathrm{U}(n)}) & (\mathbb{C}H^n,\mathrm{PU}(n,1))\\
 \mathbb H & (\mathbb{H}P^n,\mathrm{Sp}(n+1)) & (\mathbb{H}^n,\overline{\mathrm{Sp}(n) \cdot \mathrm{Sp}(1)}) & (\mathbb{H}H^n,\mathrm{Sp}(n,1))\\
 \mathbb O & (\mathbb{O}P^2,\mathrm{F}_4) & (\mathbb{O}^2,\overline{\mathrm{Spin}(9)}) & (\mathbb{O}H^2,\mathrm{F}_4^{-20})
 \end{array}
\end{equation}
Here $\mathbb{H}$ stands for the skew field of quaternions and $\mathbb{O}$ for the division algebra of octonions. In the second column stand the simply connected positively curved space forms, in the third column the corresponding flat space and in the last column the corresponding hyperbolic space form. Scaling by the curvature yields in all four cases a family $(M_\lambda,G_\lambda)$ indexed by the curvature.  

\begin{theorem}[\cite{howard93, bernig_fu_solanes}]
\label{thm_transfer}
 Let $(M_\lambda,G_\lambda), \lambda \in \R$ be one of the families from the table. Then the coalgebras $(\mathcal{C}(M_\lambda)^{G_\lambda},K_{G_\lambda}), \lambda \in \R$ are naturally isomorphic to each other. 
\end{theorem}

In particular, knowing the local kinematic formulas in the flat case, we can write down the local kinematic formulas in the curved cases, and from these formulas we can derive the global kinematic formulas in the curved cases by using \eqref{eq_local_and_global_formulas_curved}. 

\subsection{Fundamental theorem of algebraic integral geometry II:curved case}

We now formulate an analogue of Theorem \ref{thm_ftaig} in the curved case. Let $(M,G)$ be an isotropic space. Then $\mathcal{V}(M)^G$ is finite-dimensional and admits a product structure. Hence the horizontal maps in \eqref{eq_ftaig} generalize to the curved case. It is less obvious how to substitute the vertical map, i.e. the Poincar\'e duality. 

In the compact case, we can use \eqref{eq_poincare_compact} to define a map $\mathrm{PD}:\mathcal{V}(M)^G \to \mathcal{V}(M)^{G *}$. To generalize this to the non-compact case, one has to give another interpretation. It turns out that there exists a unique element $\mathrm{vol}^* \in \mathcal{V}(M)^{G *}$ with the property that $\langle \mathrm{vol}^*,\mathrm{vol}\rangle=1$ and such that $\mathrm{vol}^*$ annihilates all smooth valuations given by a pair of forms of the type $(\omega,0)$. Up to some normalization factor, $\mathrm{vol}^*$ corresponds to $\int_M$ in the compact case and to $\mathrm{PD}$ in the flat case. 

In all cases, we obtain a map (called normalized Poincar\'e duality)
\begin{displaymath}
 \mathrm{pd}:\mathcal{V}^G(M) \to \mathcal{V}^G(M)^*, \langle \mathrm{pd}(\phi),\mu\rangle:=\langle \mathrm{vol}^*,\mu \cdot \phi\rangle.
\end{displaymath}

The analogue of Theorem \ref{thm_ftaig} is the following. 
\begin{theorem}
\label{thm_ftaig_curved}
 Let $(M,G)$ be an isotropic space. 
Let $m:\mathcal{V}^G(M) \otimes \mathcal{V}^G(M) \to \mathcal{V}^G(M)$ be the restricted multiplication
map, $\mathrm{pd}:\mathcal{V}^G(M) \to \mathcal{V}^G(M)^*$ the normalized Poincar\'e duality and $k_G:\mathcal{V}^G(M)
\to \mathcal{V}^G(M) \otimes \mathcal{V}^G(M)$ the kinematic coproduct. Then the following diagram commutes:
\begin{displaymath}
\xymatrixcolsep{3pc}
\xymatrix{\mathcal{V}^G(M) \ar[d]^{\mathrm{pd}} \ar[r]^-{k_G}& \mathcal{V}^G(M) \otimes \mathcal{V}^G(M)   \ar[d]^{\mathrm{pd}
\otimes \mathrm{pd}}\\
\mathcal{V}^G(M)^* \ar[r]^-{m^*} & \mathcal{V}^G(M)^* \otimes \mathcal{V}^G(M)^*}
\end{displaymath}
\end{theorem}

As before, this can be used in two ways. Information on the algebra structure on $\mathcal{V}^G(M)$ can be translated into information on the (global) kinematic operator. In the hermitian case, we will rather go the other way and translate information from the kinematic operator into information on the algebra structure.  

The ultimate goal, however, is not to compute the global formulas only, but to compute the local kinematic formulas. One reason is that they are curvature independent (see Theorem \ref{thm_transfer}), the other is that they contain more information than the global formulas. Unfortunately, an analogue of Theorem \ref{thm_ftaig_curved} for curvature measures and local kinematic formulas is unknown (and it seems unlikely to exist). Nevertheless, there is an improvement of Theorem \ref{thm_ftaig_curved} which relates local formulas and the module structure from Theorem \ref{thm_module}. 

Let $\bar m:\mathcal{V}^G(M) \otimes \mathcal{C}^G(M) \to \mathcal{C}^G(M), \mu \otimes \Phi \mapsto \mu \cdot \Phi$ denote the module product. Alternatively, we may think of this as a map $\bar m:\mathcal{C}^G(M) \to \mathcal{C}^G(M) \otimes \mathcal{V}^G(M)^*$. 

\begin{theorem} \label{thm_semilocal}
 The following diagram commutes
 \begin{equation} \label{eq_module_semilocal}
  \xymatrix{\mathcal{C}^G(M) \ar[r]^-{\bar k_G} \ar[d]^{\mathrm{id}} & \mathcal{C}^G(M) \otimes \mathcal{V}^G(M) \ar[d]^{\mathrm{id} \otimes \mathrm{pd}} \\ \mathcal{C}^G(M) \ar[r]^-{\bar m} & \mathcal{C}^G(M) \otimes \mathcal{V}^G(M)^*}
 \end{equation}
Moreover, for $\phi \in \mathcal{V}^G(M)$ and $\Phi \in \mathcal{C}^G(M)$ we have 
\begin{displaymath}
 K_G(\phi \cdot \Phi)=(\phi \otimes \chi) \cdot K_G(\Phi)=(\chi \otimes \phi) \cdot K_G(\Phi).
\end{displaymath}

 \end{theorem}

The knowledge of the module structure will thus tell us a lot, although not everything about the local kinematic formulas. 

\section{Hermitian integral geometry: the flat case}
\label{sec_hermitian_flat}

In Section \ref{sec_kinforms} we have laid out the theoretical ground for the study of the local kinematic formulas for the spaces from Table \ref{table_isotropic_spaces}. Needless to say that to work out this program in practice requires some additional effort and ideas. In fact, only in the real and in the complex case the local kinematic formulas are known so far. 

Before going to the complex case, let us say a few words about the real case. The main point here is that the globalization map $\mathrm{glob}:\Curv^{\mathrm{SO}(n)} \to \Val^{\mathrm{SO}(n)}$ is bijective. From \eqref{eq_local_and_global_formulas_curved} we can thus read off the local kinematic formulas from the global ones, and the latter are just the well-known classical Chern-Blaschke-Santal\'o formulas. Once the local formulas are determined, one can globalize again to write down global kinematic formulas for spheres and hyperbolic spaces. 

The situation in the complex (and quaternionic and octonionic) case is different. The dimension of $\Curv^{\mathrm{U}(n)}$ is roughly twice that of $\Val^{\mathrm{U}(n)}$, hence the globalization map has a large kernel and the local formulas contain much more information than the global ones. But even to work out the global formulas is not an easy task. In this section we describe the ideas leading to a complete determination of the global kinematic formulas in the flat case $(\C^n,\mathrm{U}(n))$.

\subsection{Vector space structure}

Intrinsic volumes on a Euclidean vector space can be introduced or characterized in a number of ways, for instance averaging some functional applied to projections onto linear subspaces, averaging some functional applied to intersections with affine subspaces, or by prescribing their Klain functions. 

In the complex case, one can do analogous constructions to define some interesting valuations. However, in contrast to the real case, this gives \emph{different} bases for the space of invariant valuations. 

Alesker showed in \cite{alesker_mcullenconj01} that 
\begin{equation} \label{eq_dim_un}
 \dim \Val_k^{\mathrm{U}(n)}=\min\left\{\left\lfloor
\frac{k}{2}\right\rfloor,\left\lfloor\frac{2n-k}{2}\right\rfloor\right\}+1. 
\end{equation}

Using intersections with affine complex planes, Alesker defined 
\begin{equation*}
U_{k,p}(K) := \int_{A_\C(n,n-p)} \mu_{k-2p}(  K\cap \bar E) \,
d\bar E,
\end{equation*}
where $A_\C(n,n-p)$ is the affine complex Grassmannian, endowed with a translation invariant and $\mathrm{U}(n)$-invariant measure. 

The $U_{k,p}$, as $p$ ranges over $0,1,\ldots,\min \left\{\left\lfloor \frac k
2\right\rfloor, \left\lfloor \frac {2n-k}2 \right\rfloor\right\}$, constitute a
basis of $\Val_k^{\mathrm{U}(n)}$. 

Fu \cite{fu06} renormalized these valuations by setting 
\begin{align}
 t & := \frac{2}{\pi} \mu_1=\frac{2}{\pi} U_{1,0} \in \Val_1^{\mathrm{U}(n)}\\
s & := nU_{2,1} \in \Val_2^{\mathrm{U}(n)} \label{eq_def_s_flat}
\end{align}
which implies that 
\begin{displaymath}
s^pt^{k-2p} = \frac{(k-2p)!n!\omega_{k-2p}}{(n-p)!\pi^ {k-2p}} U_{k,p},
\end{displaymath}
where the monomial on the left hand side refers to the Alesker product.

The second basis given by Alesker uses projections onto linear complex subspaces instead of intersections:  
\begin{equation*}
C_{k,q}(K) := \int_{G_\C(n,q)} \mu_k(\pi_E(K)) \, dE,
\end{equation*}
where $G_\C(n,q)$ is the complex Grassmannian. 

As $q$ ranges over all values from $n-\min \left\{\left\lfloor \frac k
2\right\rfloor, \left\lfloor \frac {2n-k}2 \right\rfloor\right\}$ to $n$, the $C_{k,q}$ constitute a basis of
$\Val_k^{\mathrm{U}(n)}$. Up to a normalizing constant, the Alesker-Fourier transform of
$U_{k,p}$ is $C_{2n-k,n-p}$. 

A third basis, consisting of \emph{hermitian intrinsic volumes}, will be particularly useful. Recall that a real subspace $E$ of $\C^n$ is called {\it isotropic} if the restriction of the symplectic form to $E$ vanishes. Then the dimension of $E$ does not exceed $n$, and an isotropic subspace of dimension $n$ is called {\it Lagrangian}. We call $E$ of type $(k,q)$ if $E$ can be written as the orthogonal sum of a complex subspace of (complex) dimension $q$ and an isotropic
subspace of dimension $k-2q$. Then $k-q \leq n$. 

\begin{theorem}
 There is a unique valuation $\mu_{k,q} \in \Val_k^{\mathrm{U}(n)}$ whose Klain function
evaluated at a subspace of type $(k,q')$ equals $\delta_{qq'}$. Moreover, 
\begin{displaymath}
 \mathbb F \mu_{k,q}=\mu_{2n-k,n-k+q}.
\end{displaymath}
\end{theorem}

There are two more interesting bases related to hermitian intrinsic volumes. The first one is quite useful from a geometric point of view. It consists of so called \emph{Tasaki valuations}. As hermitian intrinsic volumes, they are defined via their Klain function. The orbit space $G_\C(n,k)/U(n)$ can be described in terms of $\min\left\{\left\lfloor\frac{k}{2}\right\rfloor,\left\lfloor\frac{2n-k}{2}\right\rfloor\right\}$ K\"ahler angles, and the Klain function of a Tasaki valuation is an elementary symmetric polynomial in the cosines of the K\"ahler angles. 
Comparison of the Klain functions then easily yields the explicit expression
\begin{displaymath}
\tau_{k,q} =\sum_{i=q}^{\lfloor k/2\rfloor} \binom{i}{q} \mu_{k,i}.  
\end{displaymath}

The last basis to be mentioned here is very useful for computational purposes. As is explained in \cite{schuster_lnm}, the $\mathrm{SO}(2n)$-module $\Val_k(\C^n)$ may be decomposed into a multiplicity free direct sum of irreducible representations. It is well known which of these irreducible representations contain $\mathrm{U}(n)$-invariant vectors, moreover these vectors are unique up to scale. Explicitly, this yields the valuations
\begin{equation} \label{eq_def_primitive}
\pi_{k,r} = (-1)^r(2n-4r+1)!!\sum_{i=0}^r (-1)^{i}\frac{(k-2i)!}{(2r-2i)!} \frac{(2r-2i-1)!!}{(2n-2r-2i+1)!!}\,\tau_{k,i},
\end{equation}
which form the so called \emph{primitive basis}. 

An interesting new line of research was opened by Abardia-Wannerer \cite{abardia_wannerer}, who studied versions of the classical isoperimetric inequality, but with the usual intrinsic volumes replaced by unitarily invariant valuations. In small degrees ($k \leq 3$) they study which linear combinations of hermitian intrinsic volumes satisfy an Alexandrov-Fenchel-type inequality, from which several other isoperimetric inequalities can be derived.

\subsection{Algebra structure}

The next step is to write down the algebra structure, i.e. the Alesker product. Depending on which basis we take (Akesker's valuation $U_{k,p}, C_{k,p}$, hermitian intrinsic volumes, Tasaki valuations, primitive basis), this task may be more or less difficult. The main point of entry is to determine the product of $t$ and $s$ with a hermitian intrinsic volume. We will not go into the details of these computations. A more or less complete set of formulas was written down in \cite{bernig_fu_hig}.

The computation of the algebra $\Val^{\mathrm{U}(n)}$ was obtained earlier by Fu \cite{fu06}. 
\begin{theorem}[Fu \cite{fu06}]
Let $t,s$ be variables of degree $1$ and $2$ respectively. Let $f_i$ be the component of total degree $i$ in the expansion of $\log(1+t+s)$. Then 
\begin{displaymath}
 \Val^{\mathrm{U}(n)} \cong \mathbb{R}[t,s]/(f_{n+1},f_{n+2}). 
\end{displaymath}
\end{theorem}

\subsection{Global kinematic formulas}

As stated in Theorem \ref{thm_ftaig}, the product structure may be translated into kinematic formulas. To write down explicit closed formulas valid in all dimensions is not straightforward, since we have to invert some matrices. However, in the primitive basis, the corresponding matrices are anti-diagonal and therefore easy to invert. One can then write down the principal kinematic formula in explicit form as follows.  

\begin{theorem}[Principal kinematic formula, \cite{bernig_fu_hig}] 
\label{thm_pkf_flat}
The principal kinematic formula for $(\C^n, \mathrm{U}(n))$ is given by 
\begin{displaymath}
 k_{\mathrm{U}(n)}(\chi)=\sum_{k=0}^{2n} \sum_{r=0}^{\min\left\{\lfloor \frac k 2\rfloor ,\lfloor \frac{2n-k}2\rfloor\right\}} a_{n,k,r} \pi_{k,r} \otimes{\pi_{2n-k,r}}
\end{displaymath}
with 
\begin{displaymath}
 a_{n,k,r}:=\frac {\omega_k\omega_{2n-k}}{\pi^n}  \frac{(n-r)!}{8^r(2n-4r)!} \frac{(2n-2r+1)!!}{(2n-4r+1)!!}\binom n {2r}^{-1}. 
\end{displaymath}

\end{theorem}

Together with the knowledge of the product structure, other kinematic formulas, as well as kinematic formulas in other bases, may be derived from the principal kinematic formula. It is not clear which of these formulas can be written in such a neat and closed form as the formula above. 

\section{Hermitian integral geometry: the curved case}
\label{sec_hermitian_curved}

In this section, we describe more recent work from \cite{bernig_fu_solanes}. Two of the previously open major problems in hermitian integral geometry were to compute local kinematic formulas and to compute global kinematic formulas in the curved case (say on $\mathbb{C}P^n$ and $\mathbb{C}H^n$). These two problems are interrelated by Theorem \ref{thm_transfer}: the local kinematic formula is independent of the curvature and globalizes simultaneously to the different global kinematic formulas in the space forms. Roughly speaking, the knowledge of the local kinematic formulas and the knowledge of the different global kinematic formulas is equivalent. Some partial results were known before \cite{bernig_fu_solanes}: Park established local kinematic formulas in small degrees ($n \leq 3$), while Abardia-Gallego-Solanes \cite{abardia_gallego_solanes} proved Crofton-type formulas, which are special cases of the general kinematic formulas. 

We write $\mathbb{C}P^n_\lambda$ for the complex space form with holomorphic curvature $4\lambda$. If $\lambda>0$, this is the complex projective space with an appropriate scaling of the Fubini-Study metric; if $\lambda<0$ this is complex hyperbolic space with an appropriate scaling of the Bergman metric; and if $\lambda=0$, this is just the Hermitian space $\mathbb{C}^n$. The holomorphic isometry group of $\mathbb{C}P^n_\lambda$ will be denoted by $G_\lambda$. Hence $G_\lambda \cong \mathrm{PU}(n+1)$ if $\lambda>0$; $G_\lambda \cong \mathrm{PU}(n,1)$ if $\lambda<0$ and $G_\lambda \cong \mathbb{C}^n \rtimes \mathrm{U}(n)$ if $\lambda=0$.

Before describing the solution, let us write down some intermediate steps. The first step is to write down a basis for $\mathrm{Curv}^{\mathrm{U}(n)}$. It was basically achieved by Park \cite{park02}. He determined a basis 
\begin{eqnarray} \label{eq_park_basis}
B_{k,q} & (k > 2q)\\
\Gamma_{k,q} & (n>k-q)
\end{eqnarray}
of curvature measures by writing down explicit invariant differential forms on the sphere bundle. 

It is rather easy to describe the kernel of the globalization map $\mathrm{glob}_\lambda:\mathrm{Curv}^{\mathrm{U}(n)} \cong \mathcal{C}^{G_\lambda}(\mathbb{C}P^n_\lambda) \to \mathcal{V}^{G_\lambda}(\mathbb{C}P^n_\lambda)$. More precisely, this kernel is spanned by all curvature measures of the form
\begin{displaymath}
N_{k,q}+\lambda \frac{q+1}{\pi}B_{k+2,q*1}, \quad k>2q, q>k-n,
\end{displaymath}
where $N_{k,q}:=\frac{2(n-k+q)}{2n-k}(\Gamma_{k,q}-B_{k,q})$. Note that this kernel depends on $\lambda$.  In curvature $0$, both $B_{k,q}$ and $\Gamma_{k,q}$ globalize to the hermitian intrinsic volume $\mu_{k,q}$.

The next problems are
\begin{enumerate}
 \item Compute the module structure of $\mathrm{Curv}^{\mathrm{U}(n)}$ over $\Val^{\mathrm{U}(n)}$.
 \item Compute the module structure of $\mathrm{Curv}^{\mathrm{U}(n)}$ over $\mathcal{V}^{G_\lambda}(\mathbb{C}P^n_\lambda)$. 
 \item Compute the local and semi-local formulas in curvature $\lambda$.
 \item Compute the product structure of $\mathcal{V}^{G_\lambda}(\mathbb{C}P^n_\lambda)$.
\end{enumerate}
All of these problems were solved in \cite{bernig_fu_solanes}, but not in the order given here. In fact, none of these problems could be solved independently of the others, but all had to be solved simultaneously.

\subsection{Module structure}

Even if it is not obvious how to determine the module structure of $\mathrm{Curv}^{\mathrm{U}(n)}$ over $\Val^{\mathrm{U}(n)}$, some a priori information can be obtained by geometric means.

\begin{definition}\label{def:angular cm} 
Let $V$ be a euclidean space of dimension $m$. A translation-invariant curvature measure $\Phi \in \mathrm{Curv}(V)$ is called {\bf
angular} if, for any compact convex polytope $P \subset V$,
\begin{equation}\label{eq:def angular}\Phi(P, \cdot ) = \sum_{k=0}^m \sum_{F \in \mathfrak{F}_k(P)} c_\Phi({\vec F}) \angle (F,P) \left.\mathrm{vol}_k\right|_F
\end{equation}
where $c_\Phi(\vec F)$ depends only on the $k$-plane $\vec F \in G(n,k)$ parallel to $F$ and where $\angle(F,P)$
denotes the outer angle of $F$ in $P$. The space of such translation-invariant curvature
measures will be denoted by $\mathrm{Ang}(V)$.
\end{definition}

\begin{theorem}[Angularity theorem]
\label{thm_angularity}
The product of the first intrinsic volume with an angular measure is again angular. 
\end{theorem}

Let us return to the complex case. Define curvature measures $\Delta_{kq} \in \mathrm{Curv}^{U(n)}, \ \max\{0,k-n\} \leq q \leq \frac{k}{2}<n$ by 
\begin{align}
  \Delta_{kq} :&= \frac{1}{2n-k}(2(n-k+q)\Gamma_{kq}+(k-2q)B_{kq}) \label{eq_delta_kq}\\
\Delta_{2n,n} :&= \mathrm{vol}_{2n}. \label{eq_delta_2nn}
\end{align}
Thus $\Delta_{2q,q} = \Gamma_{2q,q}$ and $\Delta_{k,k-n} = B_{k,k-n}$. Then the subspace of angular elements in $\mathrm{Curv}^{U(n)}$ is precisely the span of the $\Delta_{kq}$. This gives us some a priori information about the product of $t$ with an invariant measure, but of course does not determine it completely. 

The next piece of information concerns multiplication by $s$.

\begin{proposition} \label{prop_beta_module}
 Let $\mathrm{Beta}_k \subset \mathrm{Curv}_k^{U(n)}$ be the subspace spanned by all $B_{k,q}$. Then $s \cdot \mathrm{Beta}_k \subset \mathrm{Beta}_{k+2}$.
\end{proposition}
 
The proof of this proposition, as well as the proof of the angularity theorem, only uses elementary geometric considerations. 

Apparently, the previous proposition gives us only information about the module structure of $\mathrm{Curv}^{U(n)}$ over $\Val^{\mathrm{U}(n)}$ in the flat case. In reality, the next proposition tells us that a similar statement holds true in the curved case as well. For this, we first have to define $s$ as an element of $\mathcal{V}^{G_\lambda}(\mathbb{C}P^n_\lambda)$. Suppose first that $\lambda>0$ and fix a complex hyperplane $P_\lambda$. Then we define 
\begin{displaymath}
 s:=\frac{\lambda^{n-1}n!}{\pi^n} \int_{G_\lambda} \chi(gP_\lambda \cap \cdot) dg. 
\end{displaymath}
Regarding $\lambda$ as a parameter, one can then show that this definition extends analytically to all $\lambda \in \R$ and, in the case $\lambda=0$, coincides with the definition of $s$ from \eqref{eq_def_s_flat}. 

Recall that we have a canonical identification of the space $\mathcal{C}^{G_\lambda}(\mathbb{C}P^n_\lambda)$ with $\mathrm{Curv}^{\mathrm{U}(n)}$. 

\begin{proposition} \label{prop_mult_s_ind}
The multiplication by $s$ on $\mathcal{C}^{G_\lambda}(\mathbb{C}P^n_\lambda) \cong \mathrm{Curv}^{\mathrm{U}(n)}$ is independent of $\lambda$. 
\end{proposition}

The proof is basically an application of the transfer principle. 

Propositions \ref{prop_beta_module} and \ref{prop_mult_s_ind}, together with some elementary properties, determine completely the multiplication by $s$. 

The situation for $t$ is a bit different. First we define $t$ as an element of $\mathcal{V}^{G_\lambda}(\mathbb{C}P^n_\lambda)$. In fact, there is a unique assignment of a valuation $t$ to a Riemannian manifold with the two properties that in the Euclidean case $t=\frac{2}{\pi}\mu_1$ and that the assignment is compatible with isometric immersions. This element is the generator of the so called Lipschitz-Killing algebra \cite{alesker_survey07, bernig_fu_solanes}.

For each $k$, the valuation $t^k$ can be expressed by integration over a pair of forms $(\omega,\phi)$, both of which are defined in terms of the Riemann curvature tensor. In the special case of $\mathbb{C}P^n_\lambda$, the Riemann curvature tensor is of course well-known, and modulo some combinatorial difficulties it is straightforward to express $t^k$ as the globalization of some linear combination of basic curvature measures $B_{k,q},\Gamma_{k,q}$.  

\subsection{The isomorphism theorem}

The next aim is to determine the product structure in the curved case, i.e. the algebra structure of $\mathcal{V}^{G_\lambda}(\mathbb{C}P^n_\lambda)$.

We will need a basis of $\mathcal{V}^{G_\lambda}(\mathbb{C}P^n_\lambda)$ which will play the role of hermitian intrinsic volumes. 

\begin{definition}
For $\max\{0,k-n\}\leq q\leq \frac{k}{2}\leq n $ we set
\begin{displaymath}
\mu^\lambda_{kq}:=\sum_{i\geq 0} \frac{\lambda^i(q+i)!}{\pi^i q!} \mathrm{glob}_\lambda(\Delta_{k+2i,q+i}) \in \mathcal{V}^{G_\lambda}(\mathbb{C}P^n_\lambda),
\end{displaymath}
where $\Delta_{k+2i,q+i} \in \mathrm{Curv}^{\mathrm{U}(n)}$ was defined in \eqref{eq_delta_kq}, \eqref{eq_delta_2nn}. These valuations form a basis of $\mathcal{V}^{G_\lambda}(\mathbb{C}P^n_\lambda)$.
\end{definition}

If $k>2q$, then $\mu^\lambda_{kq}$ is just the globalization of $B_{kq}$. In the case $\lambda=0$, we obtain the usual hermitian intrinsic volumes: $\mu_{kq}^0=\mu_{kq}$. 

Using that we can write both $s^it^{k-2i}$ and $\mu^\lambda_{kq}$ as globalizations of some curvature measures, it is then possible to establish a complete dictionary between these two bases. The result is 
\begin{displaymath}
 \mu_{kq}^\lambda=(1-\lambda s)\sum_{i=q}^{\lfloor \frac{k}{2}\rfloor} (-1)^{i+q} \binom{i}{q}
\frac{\pi^k}{\omega_k(k-2i)!(2i)!}v^{\frac{k}{2}-i} u^i,
\end{displaymath}
where $v:=t^2(1-\lambda s), u:=4s-v$. Here $v^{\frac{r}{2}}$ has to be understood in the sense of power series if $r$ is odd, but only finitely many terms will be non-zero for a given dimension $n$. 

From these computations follows the next, rather surprising, theorem. 
\begin{theorem}[Isomorphism theorem] \label{thm_isomorphism}
 The algebras $\mathcal{V}^{G_\lambda}(\mathbb{C}P^n_\lambda), \lambda \in \R$ are pairwise isomorphic. More precisely, the map $s \mapsto s, t \mapsto t \sqrt{1-\lambda s}$ induces an isomorphism $\Val^{\mathrm{U}(n)} \to \mathcal{V}^{G_\lambda}(\mathbb{C}P^n_\lambda)$.
\end{theorem}

A related statement concerns the principal kinematic formula for $\mathbb{C}P^n_\lambda$.
\begin{theorem}[Global principal kinematic formula]
\label{thm_pkf_curved}
 The principal kinematic formula $k_\lambda(\chi)$, expressed in the basis $\mu_{kq}^\lambda$, is independent of $\lambda$.
\end{theorem}

Since we know the principal kinematic formula in the case $\lambda=0$ (see Theorem \ref{thm_pkf_flat}), we can therefore write down it for all $\lambda$ without any additional effort. 

\subsection{Local kinematic formulas}

By Proposition \ref{prop_mult_s_ind}, multiplication by $s$ on $\mathrm{Curv}^{\mathrm{U}(n)}$ is curvature independent. This is not true for $t$. 

However, in the flat case $\lambda=0$ one can completely determine the action of $t$ on $\mathrm{Curv}^{\mathrm{U}(n)}$ by using Theorem \ref{thm_angularity} and some elementary properties (for instance that the actions of $t$ and $s$ commute). Hence the module structure of $\mathrm{Curv}^{\mathrm{U}(n)}$ over $\Val^{\mathrm{U}(n)}$ is known. It turns out a posteriori that a version of the angularity theorem is valid in the curved case as well, but an a priori proof of this theorem seems to be missing. 

We now have enough information to find out the local kinematic formulas. More precisely, we first know that the operator $K$ is curvature independent by Theorem \ref{thm_transfer}. From the global kinematic formulas (which are known by Theorem \ref{thm_pkf_curved}), we can compute $(\mathrm{glob}_\lambda \otimes \mathrm{glob}_\lambda) \circ K$. The module structure in the flat case implies the knowledge of $(\mathrm{glob}_0 \otimes \mathrm{id}) \circ K$. These properties fix $K$. 

To write down explicit formulas, we have to introduce some more notation.

\begin{definition}
Let $N_{1,0}:=\Delta_{1,0}-B_{1,0}$. Define two natural maps from $\Val^{U(n)}$ to $\mathrm{Curv}^{U(n)}$ by
\begin{displaymath}
 \mathfrak{l}(\phi):=\phi \Delta_{0,0}, \quad \mathfrak{n}(\phi):=\phi N_{1,0}.
\end{displaymath}
\end{definition}

The importance of these maps comes from the following proposition. 

\begin{proposition} \label{prop_module_curv}
The $\Val^{\mathrm{U}(n)}$-module $\mathrm{Curv}^{\mathrm{U}(n)}$ is generated by $\Delta_{0,0}$ and $N_{1,0}$. 
\end{proposition}

Since the module structure is completely known, the proof is easy. 

Theorem \ref{thm_semilocal} and Proposition \ref{prop_module_curv} imply that the local kinematic operator $K$ is completely determined by $K(\Delta_{0,0})$ and $K(N_{1,0})$. We have already sketched above how to compute $K$. The main result is the following. 

Recall the valuations $\pi_{k,r}$ from \eqref{eq_def_primitive} and the constants $a_{n,k,r}$ from Theorem \ref{thm_pkf_flat}. Define 
\begin{align}
\rho_{kr} & := \frac{2(-1)^r (2n-4r+1)!! \pi^{k-1}}{\omega_k} \Bigg(\frac{(2r-1)!!(k+1)!}{(2n-2r+1)!!(2r)!}
\sum_{i=0}^{\lfloor\frac{k-1}{2}\rfloor} \frac{(-1)^{i+1}}{(2i+3)!(k-2i-1)!}t^{k-2i-1}u^i \nonumber \\
& \quad +\sum_{i=0}^{r-1} \frac{(-1)^i(2r-2i-3)!!}{(2n-2r-2i-1)!!(2r-2i-2)!(2i+2)!}t^{k-2i-1}u^i\Bigg) \in \Val^{\mathrm{U}(n)},
\label{eq_def_rho}
\end{align}
where $u=4s-t^2$.

\begin{theorem}[Principal local kinematic formulas, \cite{bernig_fu_solanes}] \label{thm_plkf}
 \begin{align} \label{eq_plkf}
   K(\Delta_{00}) & =\sum a_{nkr} \left[\mathfrak{l}(\pi_{kr}) \otimes \mathfrak{l}(\pi_{2n-k,r})-
\mathfrak{n}(\rho_{kr}) \otimes \mathfrak{n}(\rho_{2n-k,r})\right].\\
   K(N_{10}) & =\sum a_{nkr} \big[\mathfrak{n}(\pi_{kr}) \otimes \mathfrak{l}(\pi_{2n-k,r})+\mathfrak{l}(\pi_{kr}) \otimes \mathfrak{n}(\pi_{2n-k,r})
\nonumber \\
& \quad -\mathfrak{n}(\pi_{kr}) \otimes \mathfrak{n}(\rho_{2n-k,r}) - \mathfrak{n}(\rho_{kr}) \otimes \mathfrak{n}(\pi_{2n-k,r})\big].
\end{align}
\end{theorem}


\section{Local additive kinematic formulas for hermitian area measures and tensor valuations}
\label{sec_area_meas}

In this section, we describe a recent deep theorem by Wannerer \cite{wannerer_area_measures, wannerer_unitary_module}, which gives another type of localization of the global kinematic formulas. The theory of unitarily invariant area measures has some similarities with that of unitarily invariant curvature measures.  

\subsection{Area measures and local additive kinematic formulas}

Roughly speaking, an area measure is a valuation on a Euclidean vector space with values in the space of signed measures on the unit sphere. An example is \emph{the} area measure $S_{n-1}(K,\cdot)$ in an $n$-dimensional Euclidean space, or more generally the $k$-th area measure $S_k(K,\cdot)$, compare \cite{schneider_book14}. 

\begin{definition}
 A smooth area measure on an $n$-dimensional Euclidean vector space $V$ is a functional of the form 
\begin{displaymath}
 \Phi(K,B)=\int_{N(K) \cap \pi_2^{-1}B} \omega,
\end{displaymath}
where $\omega$ is a smooth translation invariant $(n-1)$-form on the sphere bundle $SV=V \times S^{n-1}$, $K$ is a convex body, $B \subset S^{n-1}$ a Borel subset of the unit sphere and $\pi_2:SV \to S^{n-1}$ the projection on the sphere. The space of smooth area measures is denoted by $\mathrm{Area}(V)$ or just by $\mathrm{Area}$.  
\end{definition}

Plugging $B:=S^{n-1}$ into a smooth area measure, we obtain a smooth valuation. The corresponding map is called globalization map and denoted by $\mathrm{glob}$.

As was the case for $\mathrm{Curv}$, $\mathrm{Area}$ is a module over $\Val$. However, the algebra structure on $\Val$ is not the Alesker product, but the convolution product. 
\begin{theorem}[\cite{wannerer_unitary_module}]
 There is a unique module structure of $\mathrm{Area}$ over $(\Val^\infty,*)$ such that if $\phi(K)=\mathrm{vol}(K+A)$ with $A$ a smooth convex body with strictly positive curvature, then 
\begin{displaymath}
 \phi * \Phi(K,B)=\Phi(K+A,B).
\end{displaymath}
\end{theorem}

Let $G$ be a subgroup of $\mathrm{SO}(n)$ acting transitively on the unit sphere. Then $\mathrm{Area}^G$, the space of $G$-invariant area measures, is finite-dimensional. Wannerer proved the existence of local additive kinematic formulas as follows.
\begin{theorem}[\cite{wannerer_area_measures}]
 There exists a linear map $A:\mathrm{Area}^G \to \mathrm{Area}^G \otimes \mathrm{Area}^G$ such that 
\begin{displaymath}
 A(\Phi)(K,B_1;L,B_2)=\int_G \Phi(K+gL,B_1 \cap g B_2) dg.
\end{displaymath}
It is called local additive kinematic operator. 
\end{theorem}

The local additive kinematic formula, its semi-local version $\bar a:=(\mathrm{id} \otimes \mathrm{glob}) \circ A$, and the global additive kinematic formula $a_G:\Val^G \to \Val^G \otimes \Val^G$ fit into a diagram analogous to \eqref{eq_local_and_global_formulas_curved}. Moreover, the relation between the semi-local formula and the module structure is as in \eqref{eq_module_semilocal}.

\subsection{The moment map and additive kinematic formulas for tensor valuations}

In Section \ref{sec_hermitian_curved} we have described how to obtain local kinematic formulas for curvature measures in the hermitian case. One main ingredient was the passage to complex space forms and the use of the transfer principle. For area measures, this strategy does not work, since there are no area measures on space forms. Instead, Wannerer uses the moment map, which relates local additive kinematic formulas and additive kinematic formulas for tensor valuations. 

\begin{definition}
Let $\Val^{r}:=\Val \otimes \mathrm{Sym}^r$ denote the space of tensor valuations of rank $r$. The $r$-th moment map is the map 
\begin{align*}
 M^r: \mathrm{Area} & \to \Val^r\\
\Phi & \mapsto \int_{S^{n-1}} u^r d\Phi(\cdot,u)
\end{align*}
\end{definition}

\begin{theorem}[Wannerer \cite{wannerer_area_measures}]
Let $A$ be the local additive kinematic operator, and $a^{r_1,r_2}$ the additive kinematic operator for tensor valuations (see \cite{bernig_hug_lnm, bernig_hug}).
 The following diagram commutes
\begin{displaymath}
 \xymatrix{\mathrm{Area}^G \ar[r]^-A \ar[d]^{M^{r_1+r_2}} & \mathrm{Area}^G \otimes \mathrm{Area}^G \ar[d]^{M^{r_1} \otimes M^{r_2}} \\ \Val^{r_1+r_2,G} \ar[r]^-{a^{r_1,r_2}} & \Val^{r_1,G} \otimes \Val^{r_2,G}} 
\end{displaymath}
\end{theorem}

From the theorem we can derive a strategy to compute the operator $A$, i.e. the local additive kinematic formulas: we have to choose $r_1,r_2$ in such a way that $M^{r_1}$ and $M^{r_2}$ are injective, and then $a^{r_1,r_2}$ will determine $A$.  

\subsection{Hermitian case}

Park's result \cite{park02}, which was already used in Section \ref{sec_hermitian_curved} to determine the unitarily invariant curvature measures, gives the following characterization. 
\begin{proposition}
 The area measures $B_{k,q}, \Gamma_{k,q}$ with the same restrictions on $k,q$ as in \eqref{eq_park_basis} and with $k<2n$, form a basis of the space $\mathrm{Area}^{\mathrm{U}(n)}$ of smooth $\mathrm{U}(n)$-invariant area measures.   
\end{proposition}

The module structure of $\mathrm{Area}^{\mathrm{U}(n)}$ over $(\Val^{\mathrm{U}(n)},*)$ was determined earlier by Wannerer in \cite{wannerer_unitary_module}. 

\begin{proposition}
 The module structure of $\mathrm{Area}^{\mathrm{U}(n)}$ over $(\Val^{\mathrm{U}(n)},*)$ has the following properties.
\begin{enumerate}
 \item The subspace generated by all $\Gamma_{k,q}$ is a submodule.
\item If $\hat s$ denotes the Alesker-Fourier transform of $s$, then $\hat s * B_{k,q}$ is a linear combination of $B_{k',q'}$'s. 
\end{enumerate}
\end{proposition}

Both properties follow from a careful inspection of the formula from \cite{bernig_fu06}. Since the module product is compatible with the globalization map, the first property allows to determine $\hat t * \Gamma_{k,q}$ and $\hat s * \Gamma_{k,q}$, and the second property allows to determine $\hat s * B_{k,q}$. Finally, $\hat t * B_{k,q}$ can be easily written down using some Lie derivative computation. We refer to the original paper \cite{wannerer_unitary_module} for the statement of the theorem. 

Finally, let us consider the local additive kinematic formulas for unitarily invariant area measures. In order to work out the strategy sketched in the previous section, Wannerer showed in \cite{wannerer_area_measures} that the second moment map $M^2:\mathrm{Area}^{\mathrm{U}(n)} \to \Val^{2,\mathrm{U}(n)}$ is injective. He then went on to compute the relevant parts of the additive kinematic formula $a^{2,2}$ for unitarily invariant tensor valuations. He first proved an additive version of Theorem \ref{thm_ftaig} for tensor valuations, relating $a$ and the convolution product of tensor valuations (compare also \cite{bernig_hug_lnm}). The convolution product of unitarily invariant tensor valuations of rank $2$ can be computed using a formula from \cite{bernig_fu06}. This formula is much easier to use than the corresponding formula for the product of valuations from Theorem \ref{thm_product_formula} since it does not involve any fiber integrations. 

To state his theorem, Wannerer gives a precise description of the adjoint of $A$, which is a commutative associative product on the dual space $\mathrm{Area}^{\mathrm{U}(n),*}$. It turns out that this algebra is a polynomial algebra with three generators $t,s,v$. We refer to \cite{wannerer_area_measures} for the precise statement of the theorem.   

\def\cprime{$'$}

\end{document}